\title{\bf{A new example of an algebraic surface with canonical map of degree $ 16 $} }
\author{
        NGUYEN BIN\\
 }
\date{}

\newcommand{\Addresses}{{
		\bigskip
		\footnotesize
		\noindent
		Center for Mathematical Analysis, Geometry and Dynamical Systems\\
		Departamento de Matem\'{a}tica\\
		Instituto Superior T\'{e}cnico\\
		Universidade de Lisboa\\
		Av. Rovisco Pais\\
		1049-001 Lisboa\\
		Portugal.\\
		nguyenbin@tecnico.ulisboa.pt
		
	}}

\documentclass[10pt]{article}
\usepackage{amsmath, amsthm, amssymb, mathrsfs, amscd, amsthm, amscd,amsfonts,calligra,mathrsfs,adjustbox,arydshln,lscape}

\usepackage{indentfirst,url,xypic}
\usepackage[all]{xy}
\usepackage{url}
\usepackage{tikz}
\usepackage[colorlinks,plainpages]{hyperref}
\hypersetup{
	colorlinks=true,
	linkcolor=blue,
	filecolor=magenta,      
	urlcolor=cyan,
}

\DeclareMathOperator{\degree}{deg}
\DeclareMathOperator{\image}{im}
\DeclareMathOperator{\Picard}{Pic}

\newtheorem{Theorem}{Theorem}
\newtheorem{Remark}{Remark}
\newtheorem{Proposition}{Proposition }

\newcommand{\MSC}{\textbf{Mathematics Subject Classification (2010):}}
\newcommand{\Keyword}{\textbf{Keywords:}}
\begin{document}
\maketitle
\begin{abstract}  In this note, we construct a minimal surface of general type  with geometric genus $ p_g =4 $, self-intersection of the canonical divisor $ K^2 = 32 $ and irregularity $ q = 1 $ such that its canonical map is an abelian cover of degree $ 16 $ of $\mathbb{P}^1\times  \mathbb{P}^1$. \\
	
	\noindent
	\MSC{14J29}.\\
	\Keyword{ Surfaces of general type, Canonical maps, Abelian covers}
\end{abstract}
\section{Introduction}
   Let $ X $ be a minimal smooth complex surface of general type (see \cite{MR1406314} or \cite{MR2030225}) and denote by $ \xymatrix{\varphi_{\left| K_X\right| }:X \ar@{.>}[r] & \mathbb{P}^{p_g\left( X\right)-1}} $ the canonical map of $ X $, where $ p_g\left( X\right) = \dim\left( H^0\left( X, K_{X}\right) \right) $ is the geometric genus and $ \left| K_{X} \right| $ is the canonical system of $ X $. A classical result of A. Beauville \cite[\rm Proposition 4.1 and its proof]{MR553705} says that if the image of $ \varphi_{\left | K_X\right| } $ is a surface, the degree $ d $ of the canonical map of $ X $ is less than or equal to $ 36 $ and that for large invariants the degree of the canonical map is less than or equal to $ 9 $. Later, G. Xiao improved this result by showing that if the geometric genus of $ X $ is bigger than $ 132 $, the degree of the canonical map is less than or equal to $ 8 $ (see \cite[\rm Theorem 3]{MR842626}). Only few surfaces with $ d $ greater than $ 8 $ have been known so far, such as: S. L. Tan's example \cite{MR1141782} with $ d = 9 $, U. Persson's example \cite{MR527234} with $ d = 16 $  and C. Rito's examples \cite{MR3391024}, \cite{MR3619737}, \cite{MR3663791} with $ d = 12 , 16 , 24 $. There is a recent preprint \cite{2018arXiv180711854G} of C. Gleissner, R. Pignatelli and C. Rito constructing surfaces with $ d = 24 $, $ q = 1 $ and $ d = 32 $, $ q = 0 $, where $ q = \dim\left( H^1\left( X, K_{X}\right) \right) $ is the irregularity of $ X $. \\
   
   When the canonical map has degree $ 16 $, by the Bogomolov-Miyaoka-Yau inequality, the geometric genus is less than or equal to $ 5 $ and the irregularity is at most $ 2 $. Indeed, we have 
   \begin{align*}
   16\left(p_g -2 \right)   \le 16\degree\left( \image\left( \varphi_{\left| K_X\right| }\right) \right) \le K_X^2 \le 9\chi\left( \mathcal{O}_X\right)  = 9p_g -9q +9.
   \end{align*}
   \noindent
   So $ p_g \le 5 $. In addition, since $ p_g \ge 3 $, we get $ q \le 2 $. \\

   In 1977, U. Persson gave an example with $ d = 16 $, $ K^2 = 16 $, $ p_g = 3 $ and $ q = 0 $ \cite{MR527234}. And in 2017, C. Rito constructed a surface with $ d = 16 $, $ K^2 = 16 $, $ p_g = 3 $ and $ q = 2 $ \cite{MR3619737}. As far as we know, these are all the known examples of surfaces with $ d = 16 $. In this note, we construct a minimal surface of general type with $ d = 16 $, $ K^2 = 32 $, $ p_g = 4 $ and $ q = 1 $. It is worth mentioning that the canonical map of the surface constructed by U. Persson is an abelian cover of $ \mathbb{P}^2 $.  In \cite{MR3217634}, R. Du and Y. Gao showed that if the canonical map is an abelian cover of $ \mathbb{P}^2 $ of degree $ d > 8 $, then $ d = 9 $ or $ 16 $.\\
   
   In our (very simple) construction, the canonical map is an abelian cover of $ \mathbb{P}^1 \times \mathbb{P}^1  $. The following theorem is the result of this note: 
   
   \begin{Theorem}\label{the main theorem}
   	There exists a minimal surface of general type $ X $ satisfying
   	\begin{align*}
   	K_X^2 = 32, p_g\left( X\right) = 4, q\left( X\right) = 1
   	\end{align*} 
   	such that the canonical map $ \varphi_{\left| K_X  \right|}  $ is a $ \mathbb{Z}_{2}^4- $cover of $ \mathbb{P}^1 \times \mathbb{P}^1 $.   	
   \end{Theorem}

   We construct this surface by taking a  $ \mathbb{Z}_2^4- $cover of $ \mathbb{P}^1 \times \mathbb{P}^1 $ branched in some fibres of the two rulings of   $ \mathbb{P}^1 \times \mathbb{P}^1 $.  \\
   
   \noindent
   {\bf{Notation and conventions}:} All surfaces are projective algebraic over the complex numbers. Linear equivalence of divisors is denoted by $ \equiv $. The rest of the notation is standard in algebraic geometry.

\section{$ \mathbb{Z}_{2}^4- $covers}
The construction of abelian covers was studied by R. Pardini in \cite{MR1103912}. 
For details about the building data of abelian covers  we refer the reader to Section 1 and Section 2 of R. Pardini's work (\cite{MR1103912}). 

We will denote by  $ \chi_{j_1j_2j_3j_4} $ the character of $ \mathbb{Z}_{2}^4 $ defined by
\begin{align*}
\chi_{j_1j_2j_3j_4}\left( a_1,a_2,a_3,a_4\right): =  e^{\left( \pi a_1j_1\right) i}e^{\left( \pi a_2j_2\right) i}e^{\left( \pi a_3j_3\right) i}e^{\left( \pi a_4j_4\right) i}
\end{align*}
for all $ j_1,j_2,j_3,j_4,a_1,a_2,a_3,a_4 \in \mathbb{Z}_2 $. 

We can define $ \mathbb{Z}_{2}^4- $covers as follows (see \cite[\rm Theorem 2.1]{MR1103912} and \cite[\rm Remark 1.4]{MR2956036} ):
\begin{Proposition} \label{Construction of Z16 cover}
	Given $ Y $ a smooth projective surface, let $ L_{\chi} $ be divisors of $ Y $ such that $ L_{\chi} \not\equiv \mathcal{O}_Y $ for all nontrivial characters $ \chi $ of $ \mathbb{Z}_{2}^4  $ and let $ D_{\sigma} $ be effective divisors of  $ Y $ for all $ \sigma \in \mathbb{Z}_{2}^4 \setminus \left\lbrace \left(0,0,0,0 \right)  \right\rbrace  $ such that the branch divisor $B:= \sum\limits_{\sigma \ne 0}{D_{\sigma}} $ is reduced. Then $ \left\lbrace L_{\chi}, D_{\sigma} \right\rbrace_{\chi,\sigma}$ is the building data of a  normal $ \mathbb{Z}_{2}^4-$cover $ \xymatrix{f:X \ar[r]& Y} $ if and only if
	\begin{align*}
	2L_{\chi} \equiv \sum\limits_{\chi\left( \sigma\right) = -1 }{D_{\sigma}}	
	\end{align*}
	for all nontrivial characters $ \chi $ of $ \mathbb{Z}_{2}^4  $ and for all $ \sigma \in \mathbb{Z}_{2}^4 \setminus \left\lbrace \left(0,0,0,0 \right)  \right\rbrace  $.
\end{Proposition}
\medskip

For the reader's convenience, we leave here the relations of the reduced building data of $ \mathbb{Z}_{2}^4- $ covers:
$$
\begin{adjustbox}{max width=\textwidth}
\begin{tabular}{l l r r r r r r r r r r r r r r}
$ 2L_{1000} $&$ \equiv $&$  $&$  $&$  $&$  $&$  $&$  $ &$ D_{1000} $ &$ +D_{1001} $ &$ +D_{1010} $ &$ +D_{1011} $ &$ +D_{1100} $ &$ +D_{1101} $ &$ +D_{1110} $ &$ +D_{1111} $ \\
$ 2L_{0100} $&$ \equiv $&$  $&$  $&$ D_{0100} $&$ +D_{0101} $&$ +D_{0110} $&$ +D_{0111} $ &$  $ &$  $ &$  $ &$  $ &$ +D_{1100} $ &$ +D_{1101} $ &$ +D_{1110} $ &$ +D_{1111} $ \\
$ 2L_{0010} $&$ \equiv $&$ D_{0010 } $&$ +D_{0011} $&$ $&$  $&$ +D_{0110} $&$ +D_{0111} $ &$  $ &$  $ &$ +D_{1010} $ &$ +D_{1011} $ &$  $ &$  $ &$ +D_{1110} $ &$ +D_{1111} $ \\ 
$ 2L_{0001} $&$ \equiv D_{0001 } $&$  $&$ +D_{0011} $&$  $&$ +D_{0101} $&$  $&$ +D_{0111} $ &$  $ &$ +D_{1001} $ &$  $ &$ +D_{1011} $ &$  $ &$ +D_{1101} $ &$  $ &$ +D_{1111} $ 
\end{tabular}
\end{adjustbox}
$$

By \cite[\rm Theorem 3.1]{MR1103912}  if  each $D_\sigma$ is smooth and $B $ is a normal crossings divisor, then $X$  is a smooth surface. \\

From  \cite[\rm Lemma 4.2, Proposition 4.2]{MR1103912} one has:
\begin{Proposition}\label{invariants of Z16 cover}
	If \hskip 2pt $ Y $ is a smooth surface and $ \xymatrix{f: X \ar[r]& Y} $ is a smooth $  \mathbb{Z}_{2}^4- $cover with building data $ \left\lbrace L_{\chi}, D_{\sigma} \right\rbrace_{\chi,\sigma}$, the surface $ X $ satisfies the following:
	\begin{align*}
	2K_X & \equiv f^*\left( 2K_Y + \sum\limits_{\sigma \ne 0} {D_{\sigma} } \right) \\
	K^2_X &= 4\left( 2K_Y + \sum\limits_{\sigma \ne 0} {D_{\sigma} } \right)^2 \\
	p_g\left( X \right) &=p_g\left( Y \right) +\sum\limits_{\chi \ne  \chi_{0000}  }{h^0\left( L_{\chi} + K_Y \right)} \\
	\chi\left( \mathcal{O}_X \right) &= 16\chi\left( \mathcal{O}_Y \right)  +\sum\limits_{\chi \ne \chi_{0000}  }{\frac{1}{2}L_{\chi}\left( L_{\chi}+K_Y\right)}. 
	\end{align*}
\end{Proposition}

\section{Construction}
Denote by $ F = \left\lbrace 0\right\rbrace \times \mathbb{P}^1  $ and $ G =   \mathbb{P}^1 \times \left\lbrace 0\right\rbrace$ the generators of $ \Picard\left( \mathbb{P}^1 \times \mathbb{P}^1  \right)  $. Let $ D_{1000}$, $ D_{0101} $, $ D_{0100} \in \left|  2F \right|  $, $ D_{1001} \in \left|  2G \right|  $ and $ D_{1011}$, $ D_{1010}$, $D_{0111} $, $ D_{0110} \in \left|  G \right|  $ be smooth  distinct divisors of $  \mathbb{P}^1 \times \mathbb{P}^1 $. Let $ \xymatrix{f: X \ar[r] & \mathbb{P}^1 \times \mathbb{P}^1}  $ be the $ \mathbb{Z}^4_2- $cover with the branch locus $ B = \sum\limits_{\sigma \ne 0}{D_{\sigma}} $,  where $ D_{\sigma} = 0 $ for the rest of $ D_{\sigma}$. The building data is as follows:
$$
\begin{tabular}{l r r}
$ L_{0001}\equiv $&$  F  $&$ +2G$ \\
$ L_{0010}\equiv $&$     $&$ 2G $ \\
$ L_{0100}\equiv $&$  2F $&$ +G $ \\
$ L_{1000}\equiv $&$  F  $&$ +2G$ \\
$ L_{0011}\equiv $&$  F  $&$ +2G$\\
$ L_{0101}\equiv $&$  F  $&$ +2G$\\
$ L_{0110}\equiv $&$  2F $&$ +G $\\
$ L_{0111}\equiv $&$  F  $&$ +2G$ \\
$ L_{1001}\equiv $&$  2F $&$ +G $\\
$ L_{1010}\equiv $&$  F  $&$ +2G$ \\
$ L_{1011}\equiv $&$  2F $&$ +G $\\
$ L_{1100}\equiv $&$  3F $&$ +3G$\\
$ L_{1101}\equiv $&$  2F $&$ +G $\\
$ L_{1110}\equiv $&$  3F $&$ +G $ \\
$ L_{1111}\equiv $&$  2F $&$ +G $. 
\end{tabular}
$$

Since each $ D_\sigma$ is smooth and $B$ is a normal crossings divisor,  $X$ is  smooth.  From Proposition \ref{invariants of Z16 cover}, we get	
\begin{align*}
2K_{X} & \equiv f^*\left( 2F+2G \right).
\end{align*}
\noindent
This implies that $ X $ is also minimal and of general type.  Furthermore, from Proposition \ref{invariants of Z16 cover},  $X$ has the following invariants:
\begin{align*}
K_X^2= 32, p_g\left( X\right) = 4, q\left( X\right) =1, \chi\left( \mathcal{O}_X\right) = 4.
\end{align*}

We show that the canonical map is an abelian cover of degree $ 16 $. We have the $ \mathbb{Z}_2^4- $equivariant decomposition 
\begin{align*}
H^{0}\left( X, K_X\right) = H^{0}\left( \mathbb{P}^1 \times \mathbb{P}^1, K_{\mathbb{P}^1 \times \mathbb{P}^1}\right) \oplus \bigoplus_{\chi \ne  \chi_{0000}}{H^{0}\left( \mathbb{P}^1 \times \mathbb{P}^1, K_{\mathbb{P}^1 \times \mathbb{P}^1} +L_{\chi}\right)} 
\end{align*}
\noindent
where the group $ \mathbb{Z}_2^4 $ acts on $ H^{0}\left( \mathbb{P}^1 \times \mathbb{P}^1, K_{\mathbb{P}^1 \times \mathbb{P}^1} +L_{\chi}\right) $ via the character $ \chi $ (see \cite[ \rm Proposition 4.1c]{MR1103912}). 

We consider the subgroup $ H:= \left\langle \left( 0,0,0,1\right), \left( 0,0,1,0\right), \left( 1,1,0,0\right) \right\rangle  $ of $ \mathbb{Z}_2^4 $. Let $ H^\perp $ denote the kernel of the restriction map $ \xymatrix{\left( \mathbb{Z}_2^4\right)^{*} \ar[r]&H^{*}} $, where $ \left( \mathbb{Z}_2^4\right)^{*} $ and $ H^{*} $ are the character groups of $  \mathbb{Z}_2^4 $ and $ H $, respectively. We have $ H^{\perp} = \left\langle \chi_{1100} \right\rangle  $. Because 
\begin{align*}
h^0\left( K_{\mathbb{P}^1 \times \mathbb{P}^1} +L_{ \chi} \right) = 0
\end{align*}
for all $ \chi \notin H^\perp $, the subgroup $ H $ acts trivially on $ H^{0}\left( X, K_X\right) $. So the canonical map $ \varphi_{\left| K_X \right| } $ is the composition of the quotient map $g: \xymatrix{X \ar[r]& X_1:= X/H} $ with the canonical map $\varphi_{\left| K_{X_1} \right| }$ of $ X_1 $ (see e.g. \cite[\rm Example 2.1]{MR1103913}).

The surface $ X_1 $ is the $ \mathbb{Z}_2- $cover $ \xymatrix{f_1: X_1 \ar[r] & \mathbb{P}^1 \times \mathbb{P}^1} $ branched on
\begin{align*}
B_1: &= D_{0100}+ D_{0101}+D_{0110}+D_{0111}+D_{1000}+D_{1001}+D_{1010}+D_{1011}\\
&\equiv 2L_{1100} \\
&\equiv 6F + 6G.
\end{align*}
So $ X_1  $ is a surface of general type with $ p_g\left( X_1\right) =4 $ whose only singularities are $ 36 $ nodes. Since
\begin{align*}
K_{X_1} &\equiv f_1^*\left( F +G\right), 
\end{align*}
\noindent
the linear system $ \left| K_{X_1} \right|  $ is base point free and the map $ \xymatrix{f_1: X_1 \ar[r] & \mathbb{P}^1 \times \mathbb{P}^1} $ is the canonical map of $ X_1 $. 

As the quotient map $g: \xymatrix{X \ar[r]& X_1} $ ramifies only on nodes, the canonical class  $  K_{X}   $ is the pullback of $  K_{X_1} $ and so the linear system $ \left| K_{X} \right|  $ is base point free. Therefore the canonical map $ \varphi_{\left| K_{X} \right| } $ coincides with the $ \mathbb{Z}_{2}^4- $cover $ \xymatrix{f: X \ar[r] & \mathbb{P}^1 \times \mathbb{P}^1}  $.

\begin{Remark}
	The $ \mathbb{Z}_2^3- $quotient of $ X $ by the subgroup 
	$$ \left\langle \left( 1,0,0,0\right), \left( 0,1,0,0\right), \left( 0,0,0,1\right) \right\rangle  $$
	is the product $ \mathbb{P}^1\times E $, where $ E $ is an elliptic curve. Note that also the example with $ p_g = 3 $, $ d =24 $ by C. Gleissner, R. Pignatelli and C. Rito \cite{2018arXiv180711854G} is obtained as $ \mathbb{Z}_2^3-$cover of such a product. 
	\vskip 1pt
	It is easy to check that the Albanese pencil of $ X $ has genus $ 5 $.
\end{Remark}

\section*{Acknowledgments}
The author is supported by Funda\c{c}\~{a}o para a Ci\^{e}ncia e Tecnologia (FCT), Portugal  under the framework of the program Lisbon Mathematics PhD (LisMath), Programa de Doutoramento FCT. The author is deeply indebted to Margarida Mendes Lopes for all her help and thanks Carlos Rito for many interesting conversations and suggestions. Thanks are also due to the anonymous referee for his/her thorough reading of the paper and suggestions.

\Addresses

\begin{thebibliography}{10}
	
	\bibitem{MR2956036}
	{\sc Alexeev, V., and Pardini, R.}
	\newblock Non-normal abelian covers.
	\newblock {\em Compos. Math. 148}, 4 (2012), 1051--1084.
	
	\bibitem{MR2030225}
	{\sc Barth, W.~P., Hulek, K., Peters, C. A.~M., and Van~de Ven, A.}
	\newblock {\em Compact complex surfaces}, second~ed., vol.~4 of {\em Ergebnisse
		der Mathematik und ihrer Grenzgebiete. 3. Folge. A Series of Modern Surveys
		in Mathematics [Results in Mathematics and Related Areas. 3rd Series. A
		Series of Modern Surveys in Mathematics]}.
	\newblock Springer-Verlag, Berlin, 2004.
	
	\bibitem{MR553705}
	{\sc Beauville, A.}
	\newblock L'application canonique pour les surfaces de type g\'en\'eral.
	\newblock {\em Invent. Math. 55}, 2 (1979), 121--140.
	
	\bibitem{MR1406314}
	{\sc Beauville, A.}
	\newblock {\em Complex algebraic surfaces}, second~ed., vol.~34 of {\em London
		Mathematical Society Student Texts}.
	\newblock Cambridge University Press, Cambridge, 1996.
	\newblock Translated from the 1978 French original by R. Barlow, with
	assistance from N. I. Shepherd-Barron and M. Reid.
	
	\bibitem{MR3217634}
	{\sc Du, R., and Gao, Y.}
	\newblock Canonical maps of surfaces defined by abelian covers.
	\newblock {\em Asian J. Math. 18}, 2 (2014), 219--228.
	
	\bibitem{2018arXiv180711854G}
	{\sc {Gleissner}, C., {Pignatelli}, R., and {Rito}, C.}
	\newblock {New surfaces with canonical map of high degree}.
	\newblock {\em ArXiv e-prints\/} (July 2018).
	
	\bibitem{MR1103912}
	{\sc Pardini, R.}
	\newblock Abelian covers of algebraic varieties.
	\newblock {\em J. Reine Angew. Math. 417\/} (1991), 191--213.
	
	\bibitem{MR1103913}
	{\sc Pardini, R.}
	\newblock Canonical images of surfaces.
	\newblock {\em J. Reine Angew. Math. 417\/} (1991), 215--219.
	
	\bibitem{MR527234}
	{\sc Persson, U.}
	\newblock Double coverings and surfaces of general type.
	\newblock In {\em Algebraic geometry ({P}roc. {S}ympos., {U}niv. {T}roms\o ,
		{T}roms\o , 1977)}, vol.~687 of {\em Lecture Notes in Math.} Springer,
	Berlin, 1978, pp.~168--195.
	
	\bibitem{MR3391024}
	{\sc Rito, C.}
	\newblock New canonical triple covers of surfaces.
	\newblock {\em Proc. Amer. Math. Soc. 143}, 11 (2015), 4647--4653.
	
	\bibitem{MR3663791}
	{\sc Rito, C.}
	\newblock A surface with canonical map of degree 24.
	\newblock {\em Internat. J. Math. 28}, 6 (2017), 1750041, 10.
	
	\bibitem{MR3619737}
	{\sc Rito, C.}
	\newblock A surface with {$q=2$} and canonical map of degree 16.
	\newblock {\em Michigan Math. J. 66}, 1 (2017), 99--105.
	
	\bibitem{MR1141782}
	{\sc Tan, S.~L.}
	\newblock Surfaces whose canonical maps are of odd degrees.
	\newblock {\em Math. Ann. 292}, 1 (1992), 13--29.
	
	\bibitem{MR842626}
	{\sc Xiao, G.}
	\newblock Algebraic surfaces with high canonical degree.
	\newblock {\em Math. Ann. 274}, 3 (1986), 473--483.
	
\end{thebibliography}
\end{document}